\documentclass[12pt]{article}
\usepackage{amssymb}
\usepackage{amsmath}
\usepackage{array}
\usepackage[T2A]{fontenc}
\usepackage[cp1251]{inputenc}
\renewcommand{\S}{\mathhexbox278}
 
\DeclareMathOperator{\mmod}{mod}

\DeclareMathOperator{\vep}{\varepsilon}

\DeclareMathOperator{\vf}{\varphi}

\renewcommand{\le}{\operatorname{\leqslant}}
\renewcommand{\ge}{\operatorname{\geqslant}}
\multlinegap=0em \DeclareFontFamily{T1}{msb}{}
\DeclareFontShape{T1}{msb}{m}{ol}{<5> <6> <7> <8> <9> gen * msbm
<10> <10.95> <12> <14.4> <17.28> <20.74> <24.88> msbm10}{}
\DeclareSymbolFont{AMSb}{T1}{msb}{m}{ol} \multlinegap=0em
\begin{document}

\begin{center}
{\rmfamily\bfseries\normalsize Kloosterman sums with primes to composite moduli\footnote{This work was supported by the Russian Science Foundation under grant 19\,-11-00001.}}
\end{center}

\begin{center}
{\normalsize M.A.~Korolev}
\end{center}

\vspace{0.5cm}

\fontsize{11}{12pt}\selectfont

\textbf{Abstract.} We obtain a new estimate for Kloosterman sum with primes $p\le X$ to composite modulo $q$, that is, for the exponential sum of the type
\[
\sum\limits_{p\le X,\;p\,\nmid q}\exp{\biggl(\frac{2\pi i}{q}\bigl(a\overline{p}+bp\bigr)\,\biggr)},\quad (ab,q)=1,\quad p\overline{p}\equiv 1\pmod{q},
\]
which is non-trivial in the case when $q^{\,3/4+\vep}\le X\ll q^{\,3/2}$. We also apply this estimate to the proof of solvability of some congruences with inverse prime
residues $\pmod{q}$.

\vspace{0.2cm}


\vspace{0.2cm}

\fontsize{12}{15pt}\selectfont

\begin{flushright}
\textit{To Dmitry Aleksandrovich Popov \\ on the occasion with his 80th anniversary}
\end{flushright}

\vspace{0.3cm}

\textbf{1. Introduction.}

\vspace{0.3cm}

Kloosterman sums over prime numbers has the form
\begin{equation}\label{lab_1-01}
W_{q}(a,b;X)\,=\,\mathop{{\sum}'}\limits_{p\le X}e_{q}(a\overline{p}+bp).
\end{equation}
Here $q\ge 3$, $a,b$ are integers, $(a,q)=1$, $X>1$ and $e_{q}(u) = e^{2\pi iu/q}$. The prime sign means that $p\nmid q$. For $n$ coprime to $q$, by $\overline{n} = 1/n$ we denote the inverse residue, that is, the solution of the congruence $n\overline{n}\equiv 1\pmod{q}$. Estimating (\ref{lab_1-01}), one can pursue two aims, which do not connected directly to each other, however. The first one is the estimation of $W_{q}(a,b;X)$ in the case when the length $X$ of the summation interval (as a function of $q$) is possibly smaller. The second one tries to make the decreasing factor in the estimates possibly smaller. As a rule, the results of the second type are most useful for the applications.

Thus, E.~Fouvry and P.~Michel \cite{Fouvry_Michel_1998} treat the general sum
\begin{equation}\label{lab_1-02}
W_{q}(f;X)\,=\,\mathop{{\sum}'}\limits_{p\le X}e_{q}(f(p)),
\end{equation}
where $f(x)\equiv P(x)/Q(x)$ is the rational function modulo $q$ which differs from the constant and linear function, $P, Q$ are some monic polynomials with integer coefficients.
The sum (\ref{lab_1-01}) is a partial case of such sum with $P(x)\equiv ax^{2}+b$, $Q(x)\equiv x$. In the case of prime modulo $q$, the general sum (\ref{lab_1-02}) is estimated in \cite{Fouvry_Michel_1998} as
\begin{equation}\label{lab_1-03}
W_{q}(f;X)\,\ll\,Xq^{\,\vep}\bigl(q^{\,6/7}X^{-1}\bigr)^{\! 7/32}
\end{equation}
for any $\vep>0$. This estimate is non-trivial for $X\ge q^{\,6/7+\delta}$, $\delta = \delta(\vep)>0$, and yields the asymptotic formula for the number of solutions of the congruence
\begin{equation}\label{lab_1-04}
f_{1}(p_{1})+\ldots + f_{33}(p_{33})\,\equiv\,m\pmod{q}
\end{equation}
in primes $p_{j}\le q$ for any $m$ and for any fixed tuple of the rational functions $f_{j}(x)$ satisfying the above conditions.
At the same time, in the particular case $f(x)\equiv ax^{-k}+bx$, $k\ge 1$, E.~Fouvry and P.~Michel \cite{Fouvry_Michel_1998} find the estimate
\begin{equation}\label{lab_1-05}
W_{q}(f;X)\,\ll\,Xq^{-\delta},
\end{equation}
which is valid for $q^{\,3/4+\vep}\le X\le q$. Here $\delta = \delta(\vep)$ is a sufficiently small constant (its precise dependence on $\vep$ does not specified).
The estimate (\ref{lab_1-05}) allows one to establish the solvability of the congruences similar to (\ref{lab_1-04}) in prime numbers $p_{1},\ldots, p_{k}$ lying in the ``short'' interval of the length
$X\ge q^{\,3/4+\vep}$ when $k$ is sufficiently large: $k\ge k_{0}(\vep)$.

In 2005, J.~Bourgain \cite{Bourgain_2005} proved the estimate of the sum $W_{q}(a,b;X)$ of the type (\ref{lab_1-05}) for prime $q$ and $X\ge q^{\,1/2+\vep}$. The precise formula for the power $\delta = \delta(\vep)$ (namely, $\delta = 0.0005\vep^{4}$) was given by R.C.~Baker \cite{Baker_2012}. As in the previous example, Bourgain's bound implies the solvability of the corresponding congruences in primes lying in a very short interval: $p_{1},\ldots, p_{k}\le X$, $X\ge q^{\,1/2+\vep}$. However, the number of variables should be very large: $k\gg \vep^{-4}$.

As a rule, the case of ``homogeneous'' sum
\[
W_{q}(a,0;X)\,=\,W_{q}(a,X)\,=\,\sum\limits_{p\le X}e_{q}(a\overline{p})
\]
is more easy and leads to the more precise estimates (see, for example, \cite[Appendix, Lemma 5]{Bourgain_2005}, \cite[Theorem 1]{Karatsuba_1995}, \cite[Theorem 1]{Korolev_2016}, \cite[Theorems 1, 2]{Korolev_2018a} and \cite{Korolev_2019}).

In 2010, M.Z.~Garaev \cite{Garaev_2010} proved the estimate
\begin{equation}\label{lab_1-06}
W_{q}(a;X)\,\ll\,\bigl(X^{15/16}\,+\,X^{2/3}q^{\,1/4}\bigr)q^{\vep}
\end{equation}
for prime $q$, which is non-trivial for $q^{\,3/4+\delta}\le X\le q$, $\delta = \delta(\vep)$. Using (\ref{lab_1-06}), he showed that the congruence
\[
p_{1}(p_{2}+p_{3})\,\equiv\,\lambda \pmod{q}
\]
is solvable for any $\lambda\not\equiv 0\pmod{q}$ in primes $p_{1}, p_{2}, p_{3}\le X$ for $X\ge q^{16/17+\vep}$ and thus improved the similar result of J.B.~Fridlander, P.~Kurlberg and I.E.~Shparlinski \cite{Fridlander_Kurlberg_Shparlinski_2008}. The estimate (\ref{lab_1-06}) was adopted to the case of any composite modulus $q$ and to more wide interval $q^{\,3/4+\vep}\le X\le q^{\,4/3-\vep}$ by E.~Fouvry and I.E.~Shparlinski \cite{Fouvry_Shparlinski_2011}. The last result allows one to study the divisors of the quadratic form $p_{1}p_{2}+p_{1}p_{3}+p_{2}p_{3}$ with primes lying in the interval $X<p_{1}, p_{2}, p_{3}\le 2X$, $X\to +\infty$ (see \cite{Fouvry_Michel_1998}, \cite{Baker_2012} and \cite{Korolev_2018b}).

In \cite{Korolev_2018a}, the author proved the estimate
\begin{equation}\label{lab_1-07}
W_{q}(a;X)\,\ll\,X^{32/37}q^{\,7/74+\vep}\,\ll\,X(q^{\,7/10}X^{-1})^{\,5/37}q^{\vep},
\end{equation}
which is valid for any composite modulus $q$. Thus, (\ref{lab_1-07}) expands the domain where the sum $W_{q}(a;X)$ can be estimated non-trivially, from $q^{\,3/4+\vep}$ to $q^{\,7/10+\vep}$.

In this paper, we obtain the estimate for ``non-homogeneous'' Kloosterman sum (\ref{lab_1-01}) which is valid for any composite modulus $q$, for any integers $a,b$ coprime to $q$ and for any $X$,
$q^{\,3/4+\vep}\le X\ll q^{\,3/2}$ (Theorem 1). Then we apply such estimate to the questions concerning the solvability of some congruences with prime variables lying in the ``short'' interval of the type $(1,X]$, $X\le q^{1-c}$, $c>0$ (Theorems 2, 3).

The main result is the following
\vspace{0.3cm}

\textsc{Theorem 1.} \emph{Let $0<\vep<0.1$ be any fixed constant, $q\ge q_{0}(\vep)$, $(ab,q)=1$. Then, for any $X$ satisfying the conditions $q^{\,3/4+\vep}\le X\le\,(q/2)^{3/2}$, the sum
\[
T_{q}(X)\,=\,T_{q}(a,b;X)\,=\,\mathop{{\sum}'}\limits_{n\le X}\Lambda(n)e_{q}(a\overline{n}+bn)
\]
obeys the estimate $T_{q}(X)\ll Xq^{\vep}\Delta$, where}
\begin{equation*}
\Delta\,=\,
\begin{cases}
\bigl(q^{\,3/4}X^{-1}\bigr)^{1/7}, & \textit{if}\quad q^{\,3/4}\le X\le q^{\,7/8},\\
\bigl(q^{\,2/3}X^{-1}\bigr)^{3/35}, & \textit{if}\quad q^{\,7/8}\le X\le\,(q/2)^{3/2}.
\end{cases}
\end{equation*}
\vspace{0.3cm}

\textsc{Notations.} By $\Lambda(n)$ we denote von Mangoldt function; $\tau(n)$ is the divisor function (that is, the number of divisors of $n$) and $\omega(n)$ denotes the number of different prime divisors of $n\ge 2$; for integers $a,b$ the symbol $(a,b)$ denotes the greatest common divisor of $|a|$ and $|b|$, while the symbol $(a;b)$ stands for a pair of numbers $a,b$. Prime sign in the sum means that the summation is taken over the numbers coprime to $q$.

\vspace{0.3cm}

\textbf{2. Auxiliary assertions.}

\vspace{0.3cm}

In this section, we put some auxiliary assertions. Lemmas 1, 2 are well-known. At the same time, we give here theirs proofs -- partially for the convenience of the reader, partially because we can not find this assertions in the literature in the form necessary for our purposes.

\vspace{0.3cm}

\textsc{Definition 1.} For arbitrary integers $q\ge 2$ and $A$, by $\nu(q;A)$  and $\mu(q;A)$ we denote the number of solutions of the congruences
\begin{equation}\label{lab_2-01}
x^{2}\,\equiv\,A\pmod{q}
\end{equation}
and
\begin{equation}\label{lab_2-02}
x(x+1)\,\equiv\,A\pmod{q}
\end{equation}
with the condition $1\le x\le q$, consequently.

\vspace{0.3cm}

Obviously, for a fixed $A$ these functions are multiplicative in $q$.

\vspace{0.3cm}

\textsc{Lemma 1.} \emph{Let $q = p^{\,\alpha}$, where $p\ge 2$ is prime, $\alpha\ge 1$, and let $(A,q) = p^{\,\beta}$, where $0\le \beta\le\alpha$. Then
\[
\nu(q;A)\,\le\,2^{\,c+1}p^{\,[\beta/2]},
\]
where $c = 1$ for $p = 2$ and $c = 0$ for $p\ge 3$.}

\vspace{0.3cm}

\textsc{Proof.} If $\beta = 0$  then this assertion follows from Theorems 223 and 225 of \cite{Buchstab_1966}. If $1\le \beta\le \alpha-1$, then $A = p^{\beta}a$ where $(a,p)=1$. Obviously, (\ref{lab_2-01}) is unsolvable for odd $\beta$. Hence we may assume that $\beta = 2\gamma$. Then any solution of (\ref{lab_2-01}) is divisible by $p^{\gamma}$. Setting $x = yp^{\gamma}$ we get
\[
y^{2}\equiv a\pmod{p^{\alpha-2\gamma}}.
\]
The last congruence has at most $2^{c+1}$ solutions modulo $p^{\alpha-2\gamma}$. Let $y\equiv y_{0} \pmod{p^{\alpha-2\gamma}}$ be such a solution. Then it generates the solutions of (\ref{lab_2-01}) of the type
\[
x\,\equiv\, p^{\gamma}(y_{0}+tp^{\alpha-2\gamma})\,\equiv\,p^{\gamma}y_{0}\,+\,tp^{\alpha-\gamma}\pmod{p^{\alpha}},
\]
which are not congruent $\pmod{p^{\alpha}}$ to each other for $0\le t\le p^{\gamma}-1$. Hence, $\nu(q;A)\,\le\,2^{c+1}p^{\gamma} = 2^{c+1}p^{\,[\beta/2]}$. Finally, if $\beta = \alpha$ then (\ref{lab_2-01}) is equivalent to $x^{2}\equiv 0\pmod{p^{\alpha}}$. Its solutions have the form
\[
x\,\equiv\, p^{\gamma}y\pmod{p^{\alpha}},\quad \frac{\alpha}{2}\,\le\,\gamma\,\le\,\alpha,\quad 1\le y\le p^{\alpha-\gamma},\quad (y,p)=1.
\]
Therefore, the number of such solutions is equal to
\[
\sum\limits_{\alpha/2\le\gamma\le \alpha}\varphi(p^{\alpha-\gamma})\,=\,(p-1)\sum\limits_{\alpha/2\le\gamma\le \alpha-1}p^{\alpha-\gamma-1}\,+\,1\,=\,p^{\,[\alpha/2]}\,=\,p^{\,[\beta/2]}.
\]
Then the assertion follows. $\Box$
\vspace{0.3cm}

\textsc{Corollary 1.} \emph{For any integers $q\ge 2$, $A$ one has}
\[
\nu(q;A)\,\le\,2^{\,\omega(q)+c}\sqrt{(A,q)}.
\]

\textsc{Corollary 2.} \emph{For any integers $q\ge 2$ and $A$ one has}
\[
\mu(q;a)\,\le\,4\cdot 2^{\,\omega(q)}\sqrt{(A,q)},\quad\textit{where}\quad A\,=\,4a+1.
\]
\textsc{Proof.} Multiplying both parts of (\ref{lab_2-02}) to 4 and adding 1, we obtain the congruence
\begin{equation}\label{lab_2-03}
y^{2}\,\equiv\,4a+1\pmod{4q},
\end{equation}
where $y = 2x+1$. Obviously, $\mu(q;a)$ does not exceed the number of solutions of (\ref{lab_2-03}). Therefore, by Corollary 1,
\[
\mu(q;a)\,\le\,\nu(4q;A)\,\le\,2^{\,\omega(4q)+1}\sqrt{(A,4q)}\,\le\,2^{\,\omega(q)+2}\sqrt{(A,4q)}.
\]
Since
\[
(A,4q)\,=\,(4a+1,4q)\,=\,(4a+1,q)\,=\,(A,q),
\]
we get the desired assertion. $\Box$

\vspace{0.3cm}

\textsc{Definition 2.} Let $q\ge 2$, $a,b$ be integers, and suppose that $(ab,q)=1$. Then for any $x$ coprime to $q$ we define $g(x)\equiv a\overline{x}+bx\pmod{q}$.

\vspace{0.3cm}

\textsc{Lemma 2.} \emph{Suppose that $(ab,q)=1$. Then the number $\kappa(q)$ of solutions of the congruence}
\begin{equation}\label{lab_2-04}
g(x)\,\equiv\,g(y)\pmod{q}
\end{equation}
\emph{with the conditions $1\le x,y\le q$ satisfies the estimate $\kappa(q)\le 2^{\,\omega(q)+1}\tau(q)q$.}
\vspace{0.3cm}

\textsc{Proof.} After obvious calculations we find that $\kappa(q)$ equals to the number of solutions of the congruence
\begin{equation}\label{lab_2-05}
(y-x)(y-a\overline{b}\overline{x})\,\equiv\,0\pmod{q},
\end{equation}
where $1\le x,y\le q$, $(xy,q)=1$. Since $\kappa(q)$ is a multiplicative function of $q$, it is sufficient to obtain the desired estimate in the case $q = p^{\,k}$. If $k=1$ then $\kappa(p)\le 2(p-1)$. Suppose that $k\ge 2$. Then we split all the solutions of (\ref{lab_2-05}) into the classes $E_{s}$, $0\le s\le k$. Namely, $E_{0}$ contains the solutions $(x,y)$ with the condition  $y-a\overline{b}\,\overline{x}\equiv 0 \pmod{p^{\,k}}$, $E_{k}$ contains the solutions with the condition $y-x\equiv 0\pmod{p^{\,k}}$. If $1\le s\le k-1$ then we put into the class $E_{s}$ the solutions satisfying the conditions
\begin{equation}\label{lab_2-06}
\begin{cases}
y-x\equiv 0\pmod{p^{\,s}},\\
y-a\overline{b}\overline{x}\equiv 0\pmod{p^{\,k-s}},
\end{cases}
\quad (xy,p)=1
\end{equation}
(of course, two classes can have non-empty intersection). Obviously, $|E_{0}|=\vf(p^{k})<p^{k}$, and the same is true for $|E_{k}|$. Suppose that $1\le s\le k/2$ and therefore $s\le k-s$. Then  (\ref{lab_2-06}) implies the condition
\begin{equation}\label{lab_2-07}
x^{2}\equiv a\overline{b}\pmod{p^{s}}.
\end{equation}
If $ab$ is non-residual modulo $p$ then (\ref{lab_2-07}) has no solutions. Otherwise, (\ref{lab_2-07}) has at most $e_{2}(p^{s})$ solutions modulo $p^{s}$ where $e_{n}(h)$ denotes the number of solutions of the congruence $z^{n}\equiv 1 \pmod{h}$ (see, for example, \cite{Korolev_2010}). Hence, there are at most $e_{2}(p^{s})p^{k-s}$ residues $x\pmod{p^{k}}$ satisfying (\ref{lab_2-07}). If $x_{0}$ is such a residue then (\ref{lab_2-06}) implies that
\begin{equation}\label{lab_2-08}
\begin{cases}
y\equiv x_{0}\pmod{p^{s}},\\
y\equiv a\overline{b}\overline{x}_{0}\pmod{p^{k-s}}.
\end{cases}
\end{equation}
In view of (\ref{lab_2-07}), first congruence follows from the second one. Since the second congruence has $p^{\,s}$ solutions modulo $p^{\,k}$, we get
\[
|E_{s}|\,\le\,e_{2}(p^{\,s})p^{\,k-s}\cdot p^{\,s}\,=\,e_{2}(p^{\,s})p^{\,k}.
\]
Similar arguments lead to the bound $|E_{s}|\le e_{2}(p^{\,k-s})p^{\,k}$ for $k-s<s\le k-1$. Thus, $\kappa(q)<2p^{\,k}+\sigma p^{\,k}$, where the sum $\sigma$ equals to
\[
e_{2}(p)+\ldots e_{2}(p^{m-1})+e_{2}(p^{m})+e_{2}(p^{m-1})+\ldots+e_{2}(p)\,=\,2(e_{2}(p)+\ldots + e_{2}(p^{m-1}))\,+\,e_{2}(p^{m})
\]
for even $k = 2m$ and equals to
\[
2(e_{2}(p)+\ldots + e_{2}(p^{m}))
\]
for odd $k = 2m+1$. Since
\begin{equation*}
e_{2}(p^{\,\nu})\,=\,
\begin{cases}
1,\quad & \text{if}\quad p = 2,\nu = 1,\\
2,\quad & \text{if}\quad p = 2,\nu = 2 \quad \text{or}\quad p\ge 3,\nu\ge 1,\\
4,\quad & \text{if}\quad p = 2,\nu \ge 3
\end{cases}
\end{equation*}
(see \cite[Lemma 3]{Korolev_2010}), we get
\begin{align*}
& \kappa(p^{\,k})\,<\,2\tau(p^{\,k})p^{\,k},\quad p\ge 3,\\
& \kappa(2^{\,k})\,<\,4\tau(2^{\,k})\,2^{\,k}.
\end{align*}
Thus the assertion follows. $\Box$
\vspace{0.3cm}

\textsc{Lemma 3.} \emph{Let $q\ge q_{0}, N, N_{1}$ be integers, and suppose that $N_{1}\le q$, $1<N<N_{1}\le cN$, where $c>1$ is an absolute constant. Next, denote by $I_{q}(N)$ the number of solutions of the system}
\begin{equation*}
\begin{cases}
\overline{x}_{1}+\overline{x}_{2}\equiv \overline{y}_{1}+\overline{y}_{2}\pmod{q},\\
x_{1}+x_{2}\equiv y_{1}+y_{2}\pmod{q}
\end{cases}
\end{equation*}
\emph{with the conditions $N<x_{1}, x_{2}, y_{1}, y_{2}\le N_{1}$. Then}
\[
I_{q}(N)\,<\,(2c)^{3}2^{\,\omega(q)}\tau_{3}(q)N^{2}.
\]
\textsc{Proof.} Suppose that $q$ is odd. Then we split the set of pairs $(y_{1}; y_{2})$, $N<y_{1},y_{2}\le N_{1}$, into the classes $E(\Delta,\delta)$. Namely, we put the pair $(y_{1};y_{2})$ in $E(\Delta,\delta)$ if and only if
\[
(y_{1}-y_{2},q)\,=\,\Delta,\quad (y_{1}+y_{2},q)\,=\,\delta.
\]
Obviously, the class $E(\Delta,\delta)$ is non-empty only in the case $(\Delta,\delta)=1$. Indeed, if $q$ has a prime divisor $p$ that divides both the numbers $\Delta$ and $\delta$ then
\[
y_{1}-y_{2}\,\equiv\,0\pmod{p},\quad y_{1}+y_{2}\,\equiv\,0\pmod{p}
\]
and therefore $y_{j}\equiv 0\pmod{p}$, $j = 1,2$. Since $p|q$, we have $(y_{j},q)\ne 1$, and this is impossible. Thus, $(\Delta,\delta)=1$ and $\Delta\delta|q$. Now let us fix a pair $(y_{1};y_{2})$ from the class $E(\Delta,\delta)$, where $\Delta\ne q$, and denote
\begin{equation}\label{lab_2-09}
\lambda\,\equiv\,\overline{y}_{1}+\overline{y}_{2}\pmod{q},\quad \mu\,\equiv\,y_{1}+y_{2}\pmod{q}.
\end{equation}
Thus we get the system
\begin{equation*}
\begin{cases}
\overline{x}_{1}+\overline{x}_{2}\,\equiv\,\lambda\pmod{q},\\
x_{1}+x_{2}\,\equiv\,\mu\pmod{q},
\end{cases}
\quad\text{that is,}\quad
\begin{cases}
\lambda x_{1}x_{2}\,\equiv\,\mu\pmod{q},\\
x_{1}+x_{2}\,\equiv\,\mu\pmod{q}.
\end{cases}
\end{equation*}
Hence, $x_{j}$, $j = 1,2$, are the roots of the congruence
\begin{equation}\label{lab_2-10}
\lambda x(\mu-x)\,\equiv\,\mu\pmod{q}.
\end{equation}
Setting $z = 2x-\mu$, we transform (\ref{lab_2-10}) to
\[
\lambda z^{2}\,\equiv\,\mu(\lambda\mu-4)\pmod{q}.
\]
Since $(\mu,q) = (\lambda,q) = \delta$, we have $\lambda = \delta\ell$, $\mu = \delta m$ for some $\ell$ and $m$, where
\[
\bigl(\ell, q/\delta\bigr)\,=\,\bigl(m, q/\delta\bigr)\,=\,1\quad\text{and}\quad (\ell,m)=1.
\]
Hence,
\begin{equation}\label{lab_2-11}
z^{2}\,\equiv\,A\pmod{q/\delta},
\end{equation}
where $A = m \overline{\ell} (\lambda\mu-4) \pmod{q/\delta}$. Further, $\bigl(A,q/\delta\bigr) = (\lambda\mu-4,q/\delta)$. However,
\[
\lambda\mu-4\,\equiv\,\frac{\mu^{2}}{y_{1}y_{2}}-4\,\equiv\,\frac{\mu^{2}-4y_{1}y_{2}}{y_{1}y_{2}}\,\equiv\,\frac{(y_{1}+y_{2})^{2}-4y_{1}y_{2}}{y_{1}y_{2}}\,\equiv\,\frac{(y_{1}-y_{2})^{2}}{y_{1}y_{2}}\pmod{q},
\]
so we have
\[
\bigl(\lambda\mu-4,q/\delta\bigr)\,=\,\bigl((y_{1}-y_{2})^{2},q/\delta\bigr)\,=\,\bigl(\Delta^{2},q/\delta\bigr)\,\le\,\Delta^{2}.
\]
By Corollary 1, (\ref{lab_2-11}) has at most
\[
2\cdot 2^{\,\omega(q/\delta)}\sqrt{(A,q/\delta)}\,\le\,2^{\,\omega(q)+1}\Delta.
\]
solutions in the residual system modulo $q\delta^{-1}$ and, hence, at most
\[
2^{\,\omega(q)+1}\Delta\cdot\biggl(\,\frac{2(N_{1}-N)}{q/\delta}+1\biggr)\,\le\,
2^{\,\omega(q)+1}\Delta\cdot\biggl(\,\frac{2c_{1}N\delta}{q}+1\biggr)
\]
solutions on each interval of the form $2N-\mu<z\le 2N_{1}-\mu$ (here and below $c_{1} = c-1$).

Next, since the pair $(y_{1};y_{2})$ satisfies the system
\begin{equation}\label{lab_2-12}
\begin{cases}
y_{1}-y_{2}\,\equiv\,0\pmod{\Delta},\\
y_{1}+y_{2}\,\equiv\,0\pmod{\delta},
\end{cases}
\end{equation}
by Chinese remainder theorem we have $y_{2}\equiv w\pmod{\Delta\delta}$, where $w=w(y_{1},\Delta,\delta)$ denotes some residue $\mmod \Delta\delta$. Therefore, the number of pairs
$(y_{1};y_{2})$ in the class $E(\Delta,\delta)$ does not exceed
\[
(N_{1}-N)\biggl(\frac{N_{1}-N}{\Delta\delta}+1\biggr)\,\le\,c_{1}N\biggl(\frac{c_{1}N}{\Delta\delta}+1\biggr).
\]
Hence, the contribution to $I_{q}(N)$ coming from tuples $(x_{1},y_{1},x_{2},y_{2})$ with pairs $(y_{1};y_{2})$ from $E(\Delta,\delta)$ is estimated above by
\[
c_{1}N\biggl(\frac{c_{1}N}{\Delta\delta}+1\biggr)\cdot 2^{\,\omega(q)+1}\Delta\biggl(\frac{2c_{1}N\delta}{q}+1\biggr).
\]
Further, if $(y_{1}; y_{2})$ belongs to the class $E(\Delta,\delta)$ then (\ref{lab_2-12}) implies that
\[
\Delta\,\le\,|y_{1}-y_{2}|\,\le\,c_{1}N,\quad \delta\,\le\,y_{1}+y_{2}\,\le\,2cN.
\]
Hence, the classes $E(\Delta,\delta)$ with $\Delta\ne q$ give the contribution to $I_{q}(N)$ estimated from above by
\begin{multline*}
\sum\limits_{\substack{\Delta\delta |q \\ \Delta\le c_{1}N,\;\delta\le 2cN}}2c_{1}N\biggl(\frac{c_{1}N}{\Delta\delta}+1\biggr)\,2^{\,\omega(q)}\Delta\,\biggl(2c_{1}\,\frac{N\delta}{q}+1\biggr)\,=\\
=\,2^{\,\omega(q)+1}c_{1}N\sum\limits_{\substack{\Delta\delta |q \\ \Delta\le c_{1}N,\;\delta\le 2cN}}\biggl(2c_{1}^{2}\frac{N^{2}}{q}\,+\,c_{1}\frac{N}{\delta}\,+\,2c_{1}\frac{N\Delta\delta}{q}\,+\,\Delta\biggr)\,\le\\
=\,2^{\,\omega(q)+1}c_{1}N\biggl(2c_{1}^{2}\frac{N^{2}}{q}\tau_{3}(q)\,+\,c_{1}N\,\sum\limits_{\delta|q}\frac{1}{\delta}\sum\limits_{\Delta|\tfrac{q}{\delta\mathstrut}}1\,+\,2c_{1}N\sum\limits_{n|q}\tau(n)\,\frac{n}{q}\,+\,
c_{1}N\sum\limits_{\Delta \delta|q}1 \biggr)\,\le\\
\le\,2^{\,\omega(q)+1}c_{1}N\biggl(2c_{1}^{2}\frac{N^{2}}{q}\tau_{3}(q)\,+\,c_{1}N\sum\limits_{n|q}\frac{n}{q}\,\tau(n)\,+\,2c_{1}N\sum\limits_{n|q}\frac{n}{q}\,\tau(n)\,+\,c_{1}N\tau_{3}(q)\biggr).
\end{multline*}
Since
\[
\frac{N}{q}\,<\,1,\quad \sum\limits_{n|q}\frac{n}{q}\,\tau(n)\,\le\,\sum\limits_{n|q}\tau(n)\,=\,\tau_{3}(q),
\]
the above sum does not exceed
\[
2^{\,\omega(q)+1}c_{1}N\bigl(c_{1}N\tau_{3}(q)(2c_{1}+1)\,+\,3c_{1}N\tau_{3}(q)\bigr)\,=\,4c_{1}^{2}(c+1)2^{\,\omega(q)}\tau_{3}(q)N^{2}.
\]
The contribution coming from pairs $(y_{1};y_{2})$ with the condition $\Delta=q$, that is, coming from tuples $(x_{1},y_{1},-x_{1},-y_{1})$, does not exceed $c_{1}^{2}N^{2}$. Therefore,
\[
I_{q}(N)\,\le\,4c_{1}^{2}(c+2)2^{\,\omega(q)}\tau_{3}(q)N^{2}.
\]
Now let us consider the case of even $q$. Since all the components $x_{1}, y_{1}, x_{2}, y_{2}$ are odd, then the numbers $\lambda$ and $\mu$ defined above are even. So, we have $\lambda = 2\lambda_{1}$, $\mu = 2\mu_{1}$. Further, using the same notion $E(\Delta,\delta)$ as above, for any non-empty class $E(\Delta,\delta)$ we have $(\Delta,\delta)=1$ or $(\Delta,\delta)=2$. In the first case $q = \Delta\delta u$ for some $u$, that is, $\Delta\delta|q$; in the second case $q = \tfrac{1}{2}\Delta\delta u$ and hence $\Delta\delta|2q$. If we fix some pair $(y_{1};y_{2})$ satisfying to non-empty class $E(\Delta,\delta)$, $\Delta\ne q$, and define $\lambda, \mu$ by (\ref{lab_2-09}), we see that the components $x_{1}, x_{2}$ of solution $(x_{1}, y_{1}, x_{2}, y_{2})$ satisfy the congruence (\ref{lab_2-10}) and hence the congruence
\[
\lambda_{1}x(\mu-x)\,\equiv\,\mu_{1}\pmod{\tfrac{1}{2}\,q}.
\]
Obviously,
\[
\bigl(\lambda_{1},\tfrac{1}{2}\,q\bigr)\,=\,\tfrac{1}{2}(\lambda,q)\,=\,\tfrac{1}{2}\,\delta,\quad \bigl(\mu_{1},\tfrac{1}{2}q\bigr)\,=\,\tfrac{1}{2}\,\delta,
\]
so $\lambda_{1} = \tfrac{1}{2}\delta\ell$, $\mu_{1} = \tfrac{1}{2}\,\delta m$, where $(\ell,q\delta^{-1}) = (m,q\delta^{-1}) = 1$ and $(\ell,m)=1$. Thus we have
\[
\tfrac{1}{2}\,\delta\ell x(\mu-x)\,\equiv\,\tfrac{1}{2}\delta m \pmod{\tfrac{1}{2}\,q}\quad\text{and}\quad x(\mu-x)\,\equiv\,m\overline{\ell}\pmod{q/\delta},\quad \ell\overline{\ell}\equiv 1\pmod{q/\delta}.
\]
Multiplying both parts to $(-4)$, we find
\[
(2x-\mu)^{2}\,\equiv\,m\overline{\ell}(\delta^{2}m\ell-4)\pmod{4q\delta^{-1}}\,\equiv\,4m\overline{\ell}(\lambda_{1}\mu_{1}-1)\pmod{4q\delta^{-1}}
\]
and hence
\begin{equation}\label{lab_2-13}
z^{2}\,\equiv\,A\pmod{q\delta^{-1}},
\end{equation}
where
\[
z\,=\,x-\mu_{1},\quad A\,\equiv\,m\overline{\ell}(\lambda_{1}\mu_{1}-1)\pmod{q\delta^{-1}}, \quad \ell\overline{\ell}\,\equiv\,1\pmod{q\delta^{-1}}.
\]
Obviously, $(A,q\delta^{-1})=(\lambda_{1}\mu_{1}-1,q\delta^{-1})$. Since
\[
\lambda = 2\lambda_{1}\,\equiv\,\frac{2\mu_{1}}{y_{1}y_{2}}\pmod{q},\quad\text{then}\quad \lambda_{1}\,\equiv\,\frac{\mu_{1}}{y_{1}y_{2}}\pmod{\tfrac{1}{2}\,q}
\]
and
\begin{multline*}
\biggl(A,\frac{q}{\delta}\biggr)\,=\,\biggl(\lambda_{1}\mu_{1}-1,\frac{q}{\delta}\biggr)\,=\,\biggl(\frac{\lambda\mu}{4}-1,\frac{q}{\delta}\biggr)\,=\\
=\,\biggl(\frac{(y_{1}-y_{2})^{2}}{4y_{1}y_{2}},\frac{q}{\delta}\biggr)\,=\,\biggl(\frac{(y_{1}-y_{2})^{2}}{4},\frac{q}{\delta}\bigg)\,=\,\biggl(\frac{\Delta^{2}}{4},\frac{q}{\delta}\biggr)\,\le\,\frac{\Delta^{2}}{4}.
\end{multline*}

Thus we conclude that (\ref{lab_2-13}) has at most
\[
2^{\,\omega(q/\delta)+1}\sqrt{\frac{\Delta^{2}}{4}}\,\le\,2^{\,\omega(q)}\Delta
\]
solutions modulo $q\delta^{-1}$. Therefore, the number of solutions of (\ref{lab_2-13}) lying in the interval $N-\mu_{1}<z\le N_{1}-\mu_{1}$ is bounded by
\[
\frac{c_{1}N}{q/\delta}\,+\,1\,=\,\frac{c_{1}\delta N}{q}\,+\,1.
\]
Thus, the number of tuples $(x_{1},y_{1},x_{2},y_{2})$ satisfying the initial congruence and corresponding to all classes $E(\Delta,\delta)$ with $\Delta\ne q$, is bounded by
\begin{multline*}
\sum\limits_{\substack{\Delta\delta | 2q \\ \Delta\le c_{1}N,\;\delta\le 2cN}}2c_{1}N\biggl(\frac{2c_{1}N}{\Delta\delta}+1\biggr)\,2^{\,\omega(q)}\Delta\biggl(\frac{c_{1}\delta N}{q}+1\biggr)\,\le\\
\le\,2^{\,\omega(q)+1}c_{1}N\sum\limits_{\substack{\Delta\delta | 2q \\ \Delta\le c_{1}N,\;\delta\le 2cN}}\biggl(\frac{2c_{1}^{2}N^{2}}{q}\,+\,\frac{2c_{1}N}{\delta}\,+\,\frac{c_{1}N\Delta\delta}{q}\,+\,\Delta\biggr)\,\le\\
2^{\,\omega(q)+1}c_{1}N\sum\limits_{\substack{\Delta\delta | 2q \\ \Delta\le c_{1}N,\;\delta\le 2cN}}\biggl(\frac{2c_{1}^{2}N^{2}}{q}\,+\,\frac{2c_{1}N}{\delta}\,+\,\frac{c_{1}N}{q/(\Delta\delta)}\,+\,c_{1}N\biggr)\,\le\\
\le\,2^{\,\omega(q)+1}c_{1}N\sum\limits_{\Delta\delta|2q}(2c_{1}^{2}N\,+\,4c_{1}N)\,\le\,2^{\,\omega(q)+1}\,2c_{1}^{2}(c_{1}+2)N^{2}\tau_{3}(2q)\,\le\\
\le\,8c_{1}^{2}(c_{1}+2)\,2^{\,\omega(q)}\tau_{3}N^{2}.
\end{multline*}
It remains to note that $8c_{1}^{2}(c_{1}+2) = 8(c-1)(c^{2}-1)<8c^{3}$. Lemma is proved. $\Box$
\vspace{0.3cm}

\textsc{Lemma 3.} \emph{Let $\vep>0$ be any fixed number, and suppose that $q\ge q_{0}(\vep)$ , $1<M<\tfrac{1}{2}q$. Then the number $J_{q}(M)$ of solutions of the congruence}
\[
\overline{x}_{1}+\overline{x}_{2}\,\equiv\,\overline{y}_{1}+\overline{y}_{2}\pmod{q}
\]
\emph{with the conditions $M<x_{1},y_{1},x_{2},y_{2}\le 2M$ satisfies the inequality}
\[
J_{q}(M)\,\ll_{\,\vep}\,M^{\,2+\vep}\biggl(\frac{M^{\,3/2}}{\sqrt{q}}\,+\,1\biggr).
\]

The proof of this assertion is contained in \cite{Cilleruelo_Garaev_2011} and based on the idea of D.R.~Heath-Brown \cite{Heath-Brown_1978}.
\vspace{0.3cm}

\textsc{Lemma 4.} \emph{Let $0<\vep<0.01$ be an arbitrary small but fixed number, $q\ge q_{0}(\vep)$, $a,b$ be integers, $(ab,q)=1$, $M, M_{1}, N, N_{1}$ satisfy the conditions $1<M<M_{1}$, $1<N<N_{1}$, $M,N<\tfrac{1}{2}\,q$, $M_{1}\le 2M$, $N_{1}\le 2N$. Suppose that $\{\alpha_{m}\}$, $\{\beta_{n}\}$ are any complex sequences such that $|\alpha_{m}|\le c_{1}m^{\vep}$, $|\beta_{n}|\le c_{2}n^{\vep}$ for some constants $c_{j} = c_{j}(\vep)$, $j = 1,2$ and for any $M<m\le M_{1}$, $N<n\le N_{1}$. Then the sum}
\[
C\,=\,C(M,N)\,=\,\mathop{{\sum}'}\limits_{M<m\le M_{1}}\mathop{{\sum}'}\limits_{N<n\le N_{1}}\alpha_{m}\beta_{n}e_{q}(a\overline{m}\,\overline{n}+bmn)
\]
\emph{satisfies the following inequality:}
\[
|C|\,\le\,MNq^{\,5\vep/2}\biggl(\frac{\sqrt{q}}{M\sqrt{N}}\,+\,\frac{q}{MN^{2\mathstrut}}\biggr)^{\! 1/8}.
\]

\textsc{Proof.} For any finite sequence $c = \{c_{n}\}_{n\in \mathcal{N}}$ we denote
\[
\|c\|_{r}\,=\,\biggl(\;\sum\limits_{n\in \mathcal{N}}|c_{n}|^{r}\biggr)^{\! 1/r}\quad (r\ge 1),\quad \|c\|_{\infty}\,=\,\max_{n\in \mathcal{N}}|c_{n}|.
\]
Then, by Cauchy inequality,
\begin{multline*}
|C|^{2}\,\le\,\biggl(\;\mathop{{\sum}'}\limits_{M<m\le M_{1}}|\alpha_{m}|\biggl|\mathop{{\sum}'}\limits_{N<n\le N_{1}}\beta_{n}e_{q}(a\overline{m}\,\overline{n}+bmn)\biggr|\biggr)^{2}\,\le\\
\le\,\|\alpha\|_{2}^{2}\mathop{{\sum}'}\limits_{M<m\le M_{1}}\biggl|\mathop{{\sum}'}\limits_{N<n\le N_{1}}\beta_{n}e_{q}(a\overline{m}\,\overline{n}+bmn)\biggr|^{2}\,=\\\
=\,\|\alpha\|_{2}^{2}\mathop{{\sum}'}\limits_{M<m\le M_{1}}\biggl|\mathop{{\sum}'}\limits_{N<n_{1}, n_{2}\le N_{1}}\beta_{n_{1}}\overline{\beta}_{n_{2}}e_{q}(a(\overline{n}_{1}+\overline{n}_{2})\overline{m}+b(n_{1}+n_{2})m)\biggr|\,=\\
=\,\|\alpha\|_{2}^{2}\mathop{{\sum}'}\limits_{M<m\le M_{1}}\biggl|\sum\limits_{\lambda=1}^{q}\sum\limits_{2N<\mu\le 2N_{1}}i_{1}(\lambda;\mu)e_{q}(a\lambda \overline{m}+b\mu m)\biggr|,
\end{multline*}
where $i_{1}(\lambda;\mu)$ denotes the sum of $\beta_{n_{1}}\overline{\beta}_{n_{2}}$ over the solutions of the following system:
\begin{equation}\label{lab_2-14}
\begin{cases}
\overline{n}_{1}+\overline{n}_{2}\,\equiv\,\lambda\pmod{q},\\
n_{1}+n_{2}\,\equiv\,\mu\pmod{q},
\end{cases}
\quad N<n_{1}, n_{2}\le N_{1}.
\end{equation}
Let $\theta(m)$ be an argument of the sum over $\lambda, \mu$. Then
\begin{multline*}
|C|^{2}\,\le\,\|\alpha\|_{2}^{2}\mathop{{\sum}'}\limits_{M<m\le M_{1}}e^{-i\theta(m)}\sum\limits_{\lambda=1}^{q}\sum\limits_{2N<\mu\le 2N_{1}}i_{1}(\lambda;\mu)e_{q}(a\lambda \overline{m}+b\mu m)\,\le\\
\le\,\|\alpha\|_{2}^{2}\sum\limits_{\lambda=1}^{q}\sum\limits_{2N<\mu\le 2N_{1}}|i_{1}(\lambda;\mu)|\biggl|\mathop{{\sum}'}\limits_{M<m\le M_{1}}e^{-i\theta(m)}e_{q}(a\lambda \overline{m}+b\mu m)\biggr|.
\end{multline*}
Since $|i_{1}(\lambda;\mu)|\le\|\beta\|_{\infty}^{2}\,j_{1}(\lambda;\mu)$, where $j_{1}(\lambda;\mu)$ denotes the number of solutions of (\ref{lab_2-14}), then
\[
|C|^{2}\,\le\,\|\alpha\|_{2}^{2}\|\beta\|_{\infty}^{2}\sum\limits_{\lambda=1}^{q}\sum\limits_{2N<\mu\le 2N_{1}}j_{1}(\lambda;\mu)\biggl|\mathop{{\sum}'}\limits_{M<m\le M_{1}}e^{-i\theta(m)}e_{q}(a\lambda \overline{m}+b\mu m)\biggr|.
\]
Using Cauchy inequality again, we get
\begin{multline}\label{lab_2-15}
|C|^{4}\,\le\,\|\alpha\|_{2}^{4}\|\beta\|_{\infty}^{4}\biggl(\,\sum\limits_{\lambda=1}^{q}\sum\limits_{2N<\mu\le 2N_{1}}j_{1}(\lambda;\mu)\biggr)\times\\
\sum\limits_{\lambda=1}^{q}\sum\limits_{2N<\mu\le 2N_{1}}j_{1}(\lambda;\mu)\biggl|\mathop{{\sum}'}\limits_{M<m\le M_{1}}e^{-i\theta(m)}e_{q}(a\lambda \overline{m}+b\mu m)\biggr|^{2}.
\end{multline}
The first double sum in (\ref{lab_2-15}) coincides with the number of all pairs $(n_{1};n_{2})$ and does not exceed $N^{2}$. Therefore,
\[
|C|^{4}\,\le\,\|\alpha\|_{2}^{4}\|\beta\|_{\infty}^{4}N^{2}
\sum\limits_{\lambda=1}^{q}\sum\limits_{2N<\mu\le 2N_{1}}j_{1}(\lambda;\mu)\biggl|\mathop{{\sum}'}\limits_{M<m\le M_{1}}e^{-i\theta(m)}e_{q}(a\lambda \overline{m}+b\mu m)\biggr|^{2}.
\]
Third application of Cauchy inequality yields:
\begin{multline}\label{lab_2-16}
|C|^{8}\,\le\,\|\alpha\|_{2}^{8}\|\beta\|_{\infty}^{8}N^{4}\biggl(\sum\limits_{\lambda=1}^{q}\sum\limits_{2N<\mu\le 2N_{1}}j_{1}^{2}(\lambda;\mu)\biggr)\times\\
\sum\limits_{\lambda=1}^{q}\sum\limits_{2N<\mu\le 2N_{1}}\biggl|\mathop{{\sum}'}\limits_{M<m\le M_{1}}e^{-i\theta(m)}e_{q}(a\lambda \overline{m}+b\mu m)\biggr|^{4}.
\end{multline}
First double sum in (\ref{lab_2-16}) coincides with the number $I_{q}(N)$ from Lemma 3 (with $c = 2$). Thus we have
\begin{multline*}
|C|^{8}\,\le\,\|\alpha\|_{2}^{8}\|\beta\|_{\infty}^{8}N^{4}I_{q}(N)
\sum\limits_{\lambda=1}^{q}\sum\limits_{2N<\mu\le 2N_{1}}\mathop{{\sum}'}\limits_{M<m_{1},m_{2},m_{3},m_{4}\le M_{1}}e^{-i(\theta(m_{1})+\theta(m_{2})-\theta(m_{3})-\theta(m_{4}))}\times\\
\times e_{q}(a\lambda (\overline{m}_{1}+\overline{m}_{2}-\overline{m}_{3}-\overline{m}_{4})+b\mu (m_{1}+m_{2}-m_{3}-m_{4}))\,=\\
=\,\|\alpha\|_{2}^{8}\|\beta\|_{\infty}^{8}N^{4}I_{q}(N)
\sum\limits_{\lambda=1}^{q}\;\sum\limits_{2N<\mu\le 2N_{1}}\;\sum\limits_{\sigma=1}^{q}
\sum\limits_{|\tau|\le 2M}i_{2}(\sigma;\tau)e_{q}(a\lambda\sigma+b\mu\tau),
\end{multline*}
where $i_{2}(\sigma;\tau)$ denotes the sum of the terms $e^{-i(\theta(m_{1})+\theta(m_{2})-\theta(m_{3})-\theta(m_{4}))}$ over the solutions of the system
\begin{equation}\label{lab_2-17}
\begin{cases}
\overline{m}_{1}+\overline{m}_{2}\,\equiv\,\overline{m}_{3}+\overline{m}_{4}+\sigma\pmod{q},\\
m_{1}+m_{2}\,\equiv\,m_{3}+m_{4}+\tau\pmod{q},
\end{cases}
\quad M<m_{1}, m_{2}, m_{3}, m_{4}\le M_{1}.
\end{equation}
Setting $j_{2}(\sigma;\tau)$ for the number of solutions of (\ref{lab_2-17}), we get
\begin{multline*}
|C|^{8}\,\le\,\|\alpha\|_{2}^{8}\|\beta\|_{\infty}^{8}N^{4}I_{q}(N)\times\\
\sum\limits_{\sigma=1}^{q}\sum\limits_{|\tau|\le 2M}j_{2}(\sigma;\tau)
\biggl|\sum\limits_{\lambda = 1}^{q}\sum\limits_{2N<\mu\le 2N_{1}}e_{q}(a\lambda\sigma+b\mu\tau)\biggr|\,\le\\
\le\,\|\alpha\|_{2}^{8}\|\beta\|_{\infty}^{8}N^{4}I_{q}(N)\sum\limits_{\sigma=1}^{q}
\sum\limits_{|\tau|\le 2M}j_{2}(\sigma;\tau)\sum\limits_{2N<\mu\le 2N_{1}}\biggl|\,\sum\limits_{\lambda=1}^{q}e_{q}(a\lambda\sigma)\biggr|\,\le\\
\le\,2\|\alpha\|_{2}^{8}\|\beta\|_{\infty}^{8}N^{5}I_{q}(N)\sum\limits_{\sigma=1}^{q}\biggl(\,\sum\limits_{|\tau|\le 2M}j_{2}(\sigma;\tau)\biggr)\biggl|\,\sum\limits_{\lambda=1}^{q}e_{q}(a\lambda\sigma)\biggr|.
\end{multline*}
Since $(a,q)=1$, the inner sum is equal to $q$ if $\sigma = q$ and equals to $0$ otherwise. Next,
\[
\sum\limits_{|\tau|\le M}j_{2}(q;\tau)\,\le\,J_{q}(M).
\]
Therefore,
\[
|C|^{8}\,\le\,2\|\alpha\|_{2}^{8}\|\beta\|_{\infty}^{8}N^{5}I_{q}(N)J_{q}(M).
\]
Using the estimates of lemmas 2, 3, we get
\begin{multline*}
|C|^{8}\,\le\,2(c_{1}c_{2})^{8}M^{4+8\vep}N^{5+8\vep}q\,2^{\,\omega(q)}\tau_{3}(q)(12N)^{2}q^{\,\vep}M^{2}
\biggl(\,\frac{M^{\,3/2}}{\sqrt{q}}\,+\,1\biggr)\,\ll\\
\ll\,q^{\,10\vep}M^{6}N^{7}q\biggl(\,\frac{M^{\,3/2}}{\sqrt{q}}\,+\,1\biggr)\,\ll\,q^{\,10\vep}(MN)^{8}\,\frac{q}{M^{2}N}
\biggl(\,\frac{M^{\,3/2}}{\sqrt{q}}\,+\,1\biggr)\,\ll\\
\ll\,q^{\,10\vep}(MN)^{8}\,\biggl(\,\frac{\sqrt{q}}{N\sqrt{M}}\,+\,\frac{q}{M^{2}N}\biggr).
\end{multline*}
Lemma is proved. $\Box$
\vspace{0.3cm}

\textsc{Corollary.} \emph{Under the assumptions of lemma, the sum}
\[
S\,=\,S(M,N)\,=\,\mathop{{\sum}'}\limits_{M<m\le M_{1}}\mathop{{\sum}'}\limits_{\substack{N<n\le N_{1} \\ mn\le X}}\alpha_{m}\beta_{n}e_{q}(a\overline{m} \overline{n}+bmn),\quad X>1,
\]
\emph{satisfies the inequality}
\[
|S|\,\le\,MNq^{\,3\vep/2}\biggl(\,\frac{\sqrt{q}}{N\sqrt{M}}\,+\,\frac{q}{M^{2}N}\biggr)^{\! 1/8}.
\]
\textsc{Proof.} If $M_{1}N_{1}\le X$ then the sum $S(M,N)$ coincides the sum $C(M,N)$; if $X<MN$ then this sum is empty. Thus, we may assume that $MN\le X\le M_{1}N_{1}$. For fixed $m$, $M<m\le M_{1}$, we set $N_{2} = \min{\bigl(Xm^{-1},N_{1}\bigr)}$. Using the relation
\begin{equation*}
\frac{1}{q}\sum\limits_{|c|\le q/2}\,\sum\limits_{N<\nu\le N_{2}}e_{q}(c(n-\nu))\,=
\begin{cases}
1, & \text{if}\;\;N<n\le N_{2},\\
0, & \text{otherwise},
\end{cases}
\end{equation*}
we get
\begin{multline*}
S(M,N)\,=\,\mathop{{\sum}'}\limits_{M<m\le M_{1}}\mathop{{\sum}'}\limits_{N<n\le N_{1}}\frac{1}{q}\sum\limits_{|c|\le q/2}\;\sum\limits_{N<\nu\le N_{2}}e_{q}(c(n-\nu))\alpha_{m}\beta_{n}e_{q}(a\overline{m} \overline{n}+bmn)\,=\\
=\,\sum\limits_{|c|\le q/2}\;\mathop{{\sum}'}\limits_{M<m\le M_{1}}\sum\limits_{N<\nu\le N_{1}}
\frac{\alpha_{m}}{q}\,\biggl(\;\sum\limits_{N<\nu\le N_{2}}e_{q}(-c\nu)\biggr)
\sum\limits_{N<\nu\le N_{1}}\beta_{n}e_{q}(cn)e_{q}(a\overline{m} \overline{n}+bmn)\,=\\
=\,\sum\limits_{|c|\le q/2}\frac{S_{c}(M,N)}{|c|+1},
\end{multline*}
where
\begin{multline*}
S_{c}(M,N)\,=\,\mathop{{\sum}'}\limits_{M<m\le M_{1}}\mathop{{\sum}'}\limits_{N<n\le N_{1}}a_{m}b_{n}e_{q}(a\overline{m} \overline{n}+bmn),\\
a_{m}\,=\,\alpha_{m}\,\frac{|c|+1}{q}\sum\limits_{N<\nu\le N_{2}}e_{q}(-c\nu),\quad b_{n}\,=\,\beta_{n}e_{q}(cn).
\end{multline*}
One can check that $|a_{m}|\le |\alpha_{m}|$, $|b_{n}|=|\beta_{n}|$ for any $m,n$. Using the estimate of Lemma 4, we find that
\[
|S|\,\le\,MNq^{\,5\vep/4}\Delta\sum\limits_{|c|\le q/2}\frac{1}{|c|+1}\,\le\,MNq^{\,5\vep/4}(\ln{q}+1)\Delta\,\le\,MNq^{\,3\vep/2}\Delta,
\]
where
\[
\Delta\,=\,\biggl(\,\frac{\sqrt{q}}{M\sqrt{N}}\,+\,\frac{q}{MN^{2}}\biggr)^{\! 1/8}.\quad \Box
\]
\vspace{0.2cm}

\textsc{Lemma 5.}  \emph{Let $q\ge 2$, $a, b$ be integers, and suppose that $a$ or $b$ are not divisible by $q$. Then the following estimate holds:}
\[
\biggl|\,\mathop{{\sum}'}\limits_{n=1}^{q}e_{q}(a\overline{n}+bn)\biggr|\,\le\,\tau(q)\sqrt{q}\,(a,b,q)^{1/2}.
\]

For the proof, see \cite{Esterman_1961}.
\vspace{0.2cm}

\textsc{Corollary.} \emph{Under the assumptions of the lemma, for any $N\le q$ we have}
\[
\biggl|\,\mathop{{\sum}'}\limits_{n=1}^{N}e_{q}(a\overline{n}+bn)\biggr|\,\le\,\tau(q)\sqrt{q}\,(a,q)^{1/2}(\ln{q}+1).
\]
\vspace{0.2cm}

\textbf{3. Main theorem.}

\vspace{0.3cm}

Let $V$ satisfies the conditions $1<V<\sqrt{X}$ and $XV^{-1}\le \tfrac{1}{2}\,q$. Applying Vaughan's identity in the form given in \cite[Ch. II, \S 6, Theorem 1]{Karatsuba_Voronin_1992}, we get
\begin{multline*}
T_{q}(X)\,=\,\mathop{{\sum}'}\limits_{m\le V}\mu(m)\mathop{{\sum}'}\limits_{n\le Xm^{-1}}(\ln{n})e_{q}(a\overline{m}\overline{n}+bmn)\,-\\
-\mathop{{\sum}'}\limits_{k,\ell\le V}\mu(k)\Lambda(\ell)\mathop{{\sum}'}\limits_{n\le X(k\ell)^{-1}}e_{q}(a\overline{k}\overline{\ell} \overline{n}+bk\ell n)\,-\\
-\,\mathop{{\sum}'}\limits_{V<m\le XV^{-1}}b_{m}\mathop{{\sum}'}\limits_{V<n\le Xm^{-1}}\Lambda(n)e_{q}(a\overline{m}\overline{n}+bmn)\,+\,O(V),\quad b_{m}\,=\,\sum\limits_{d|m,\;\;d\le V}\mu(d).
\end{multline*}
Setting
\[
a_{m}\,=\,\sum\limits_{k\ell = m,\;k,\ell\le V}\mu(k)\Lambda(\ell),
\]
we find
\begin{multline*}
T_{q}(X)\,=\,\mathop{{\sum}'}\limits_{m\le V}\mu(m)\mathop{{\sum}'}\limits_{n\le Xm^{-1}}(\ln{n})e_{q}(a\overline{m}\overline{n}+bmn)\,-\\
-\mathop{{\sum}'}\limits_{m\le V}a_{m}\mathop{{\sum}'}\limits_{n\le Xm^{-1}}e_{q}(a\overline{m}\overline{n}+bmn)\,-\,\mathop{{\sum}'}\limits_{V<m\le V^{2}}a_{m}\mathop{{\sum}'}\limits_{n\le Xm^{-1}}e_{q}(a\overline{m}\overline{n}+bmn)-\,\\
-\,\mathop{{\sum}'}\limits_{V<m\le XV^{-1}}b_{m}\mathop{{\sum}'}\limits_{V<n\le Xm^{-1}}\Lambda(n)e_{q}(a\overline{m}\overline{n}+bmn)\,+\,O(V)\,=\,S_{1}-S_{2}-S_{3}-S_{4}+O(V),
\end{multline*}
where the notations $S_{j}$ are obvious.  Corollary of Lemma 5 implies
\begin{equation}\label{lab_3-01}
|S_{2}|\,\le\,\sum\limits_{m\le V}(\ln{m})\tau(q)\sqrt{q}(\ln{q}+1)\,\le\,V\sqrt{q}\,q^{\,\vep}.
\end{equation}
By Abel summation, we similarly get the estimate
\begin{equation}\label{lab_3-02}
|S_{1}|\,\le\,V\sqrt{q}\,q^{\,\vep}.
\end{equation}
Next, we split the sums $S_{3}, S_{4}$ into double sums of the type
\[
S(M,N)\,=\,\mathop{{\sum}'}\limits_{M<m\le M_{1}}\,\mathop{{\sum}'}\limits_{\substack{N<n\le N_{1} \\ mn\le X}}\alpha_{m}\beta_{n}e_{q}(a\overline{m}\overline{n}+bmn),
\]
where
\[
|\alpha_{m}|\,\le\,\ln{m},\quad |\beta_{n}|\le 1,\quad V\le M<V^{2},\quad 1\le N\le XM^{-1}
\]
in the case of the sum $S_{3}$ and
\[
|\alpha_{m}|\,\le\,\tau(m),\quad |\beta_{n}|\le \Lambda(n),\quad V\le M<XV^{-1},\quad V\le N\le XM^{-1}
\]
for the sum $S_{4}$. Further, let $1<D<X$ be some parameter to be chosen later. All the sums $S(M,N)$ with $MN\le D$ we estimate trivially:
\begin{equation}\label{lab_3-03}
|S(M,N)|\,\le\,\sum\limits_{M<m\le M_{1}}|\alpha_{m}|\;\sum\limits_{N<n\le N_{1}}|\beta_{n}|\,\ll\,MN(\ln{q})\,\ll\,Dq^{\,\vep/2}.
\end{equation}
In the case $D>MN$ we use the corollary of Lemma 4:
\begin{equation}\label{lab_3-04}
|S(M,N)|\,\le\,MNq^{\,3\vep/2}\biggl(\frac{\sqrt{q}}{N\sqrt{M}}\,+\,\frac{q}{NM^{2\mathstrut}}\biggr)^{\! 1/8}.
\end{equation}
Changing the order of summation over $m$ and $n$, we also get the inequality
\begin{equation}\label{lab_3-05}
|S(M,N)|\,\le\,MNq^{\,3\vep/2}\biggl(\frac{\sqrt{q}}{M\sqrt{N}}\,+\,\frac{q}{MN^{2\mathstrut}}\biggr)^{\! 1/8}.
\end{equation}
Now we estimate the contribution of (\ref{lab_3-04}) and (\ref{lab_3-05}) to $S_{3}$. Suppose that $M<N$. Then
\[
N^{2}>MN>D,\quad N>\sqrt{D}, \quad\text{and hence}\quad N\sqrt{M}\,=\,\sqrt{N}\sqrt{MN}\,>\,D^{\,1/4}\sqrt{D}\,=\,D^{\,3/4}.
\]
Further, $NM^{2}\,=\,MN\cdot M\,>\,DV$. Thus, (\ref{lab_3-04}) implies the inequality
\begin{equation}\label{lab_3-06}
|S(M,N)|\,\le\,MN q^{\,3\vep/2}\biggl(\frac{\sqrt{q}}{D^{3/4}}\,+\,\frac{q}{DV}\biggr)^{\! 1/8}.
\end{equation}
Suppose now that $M\ge N$. Then
\[
M^{2}\,\ge\,MN\,>\,D,\quad M\,>\,\sqrt{D},\quad \text{and hence}\quad M\sqrt{N}\,>\,D^{3/4}.
\]
Moreover,
\[
MN^{2}\,=\,\frac{(MN)^{2}}{M}\,>\,\frac{D^{2}}{M}\,>\,\frac{D^{2}}{V^{2}}.
\]
Therefore, the bound (\ref{lab_3-05}) implies the estimate
\begin{equation}\label{lab_3-07}
|S(M,N)|\,\le\,MN q^{\,3\vep/2}\biggl(\frac{\sqrt{q}}{D^{3/4}}\,+\,\frac{qV^{2}}{D^{2}}\biggr)^{\! 1/8}.
\end{equation}
Using both (\ref{lab_3-04}) and (\ref{lab_3-05}), for any $M,N$ under considering we get
\begin{equation}\label{lab_3-08}
|S(M,N)|\,\le\,MN q^{\,3\vep/2}\biggl(\frac{\sqrt{q}}{D^{3/4}}\,+\,\frac{q}{DV}\,+\,\frac{qV^{2}}{D^{2}}\biggr)^{\! 1/8}.
\end{equation}
The summation over $M,N$ together with (\ref{lab_3-03}) give
\begin{equation}\label{lab_3-09}
|S_{3}|\,\le\,X q^{\,2\vep}\biggl\{\frac{D}{X}\,+\,\biggl(\frac{\sqrt{q}}{D^{3/4}}\,+\,\frac{q}{DV}\,+\,\frac{qV^{2}}{D^{2}}\biggr)^{\! 1/8}\,\biggr\}.
\end{equation}
The sum $S_{4}$ is treated in a similar way. Since $M,N$ satisfy symmetric conditions $V<M,N\le XV^{-1}$, it is sufficient to consider only the case $M\le N$. Thus we get
\[
N\,\ge\,\sqrt{D},\quad N\sqrt{M}\,>\,D^{3/4},\quad  NM^{2}\,>\,DM\,>\,DV
\]
and therefore
\[
|S(M,N)|\,<\,MNq^{\,3\vep/2}\biggl(\frac{\sqrt{q}}{D^{3/4}}\,+\,\frac{q}{DV}\biggr)^{\! 1/8}.
\]
The summation over $M,N$ together with (\ref{lab_3-03}) imply
\begin{equation}\label{lab_3-10}
|S_{4}|\,\le\,X q^{\,2\vep}\biggl\{\frac{D}{X}\,+\,\biggl(\frac{\sqrt{q}}{D^{3/4}}\,+\,\frac{q}{DV}\biggr)^{\! 1/8}\,\biggr\}.
\end{equation}
Summing the estimates (\ref{lab_3-01}), (\ref{lab_3-02}), (\ref{lab_3-09}) and (\ref{lab_3-10}), we find
\[
T_{q}(X)\,\ll\,Xq^{2\vep}\biggl\{\frac{D}{X}\,+\,\frac{V\sqrt{q}}{X}\,+\,\biggl(\frac{\sqrt{q}}{D^{3/4}}\,+\,\frac{q}{DV}\,+\,\frac{qV^{2}}{D^{2}}\biggr)^{\! 1/8}\,\biggr\}\,\ll\, Xq^{2\vep}\Delta,
\]
where
\[
\Delta\,=\,\delta^{1/8},\quad \delta\,=\,\frac{D^{8}}{X^{8}}\,+\,\frac{V^{8}q^{4}}{X^{8}}\,+\,\frac{\sqrt{q}}{D^{3/4}}\,+\,\frac{q}{DV}\,+\,\frac{qV^{2}}{D^{2}}.
\]
Now we choose $V$ from the equation
\[
\frac{q}{DV}\,=\,\frac{qV^{2}}{D^{2}},\quad\text{that is,}\quad V = D^{1/3}
\]
(so, the condition $V<\sqrt{X}$ holds automatically). Then
\[
\delta\,\ll\,\frac{D^{8}}{X^{8}}\,+\,\frac{q^{4}D^{8/3}}{X^{8}}\,+\,\frac{\sqrt{q}}{D^{3/4}}\,+\,\frac{q}{D^{4/3}}.
\]
If we define $D$ by the relation
\begin{equation}\label{lab_3-11}
\frac{D^{8}}{X^{8}}\,=\,\frac{\sqrt{q}}{D^{3/4}},\quad\text{that is,}\quad D\,=\,q^{\,2/35}X^{\,32/35},
\end{equation}
then we find
\[
\delta\,\ll\,q^{\,16/35}X^{\,-24/35}\,+\,q^{\,97/105}X^{\,-128/105},
\]
or, that is the same,
\begin{equation}\label{lab_3-12}
\delta\,\ll\,
\begin{cases}
q^{\,97/105}X^{\,-128/105}, & \text{if}\quad q^{\,97/128}\le X\le q^{\,7/8},\\
q^{\,16/35}X^{\,-24/35}, & \text{if}\quad X\ge q^{\,7/8}.
\end{cases}
\end{equation}

An upper bound for $X$ in the last estimate in (\ref{lab_3-12}) is calculated as follows. The quantities $M,N$ in Lemma 4 obey the conditions $M,N\le \tfrac{1}{2}\,q$. At the same time, upper bounds for $M,N$ in $S(M,N)$ are $XV^{-1}$ and $V^{2}$. Hence, it is necessary to satisfy the conditions $XV^{-1}\le \tfrac{1}{2}\,q$ and $V^{2}\le \tfrac{1}{2}\,q$. Since $V = D^{1/3}$ and $D<X$, then $V^{2} = DV^{-1}<XV^{-1}$. So, it is sufficient to check the condition $XV^{-1}\le \tfrac{1}{2}\,q$ or, that is the same, $2X\le qD^{1/3}$. In view of (\ref{lab_3-11}), we get
\begin{equation}\label{lab_3-13}
X\,\le\,c_{1}q^{\,107/73}, \quad c_{1}=2^{-105/73}.
\end{equation}
By (\ref{lab_3-12}), (\ref{lab_3-13}) we get:
\begin{equation}\label{lab_3-14}
\delta\,\ll\,
\begin{cases}
q^{\,97/105}X^{-128/105},& \text{if}\quad q^{\,97/128}\,\le\,X\,\le\,q^{\,7/8},\\
q^{\,16/35}X^{-24/35},& \text{if}\quad q^{\,7/8}\,\le\,X\,\le\,c_{1}q^{\,107/73}.
\end{cases}
\end{equation}

Now let us choose $D$ from the relation
\begin{equation}\label{lab_3-15}
\frac{D^{8}}{X^{8}}\,=\,\frac{q}{D^{4/3}},\quad\text{that is,}\quad D\,=\,q^{\,3/28}X^{6/7}.
\end{equation}
Then
\[
\delta\,\ll\,q^{\,6/7}X^{-8/7}\,+\,q^{\,47/112}X^{-9/14},
\]
or, that is the same,
\begin{equation}\label{lab_3-16}
\delta\,\ll\,
\begin{cases}
q^{\,6/7}X^{\,-8/7}, & \text{if}\quad q^{\,3/4}\le X\le q^{\,7/8},\\
q^{\,47/112}X^{\,-9/14}, & \text{if}\quad X\ge q^{\,7/8}.
\end{cases}
\end{equation}
An upper bound for $X$ in the last estimate in (\ref{lab_3-16}) is defined from the condition $X\le\tfrac{1}{2}\,qV = \tfrac{1}{2}\,qD^{1/3}$. This inequality together with (\ref{lab_3-15}) imply
\begin{equation}\label{lab_3-17}
X\,\le\,c_{2}q^{\,29/20},\quad c_{2}=2^{-7/5}.
\end{equation}
Thus, from (\ref{lab_3-16}) and (\ref{lab_3-17}) we conclude that
\begin{equation}\label{lab_3-18}
\delta\,\ll\,
\begin{cases}
q^{\,6/7}X^{-8/7},& \text{if}\quad q^{\,3/4}\,\le\,X\,\le\,q^{\,7/8},\\
q^{\,47/112}X^{-9/14},& \text{if}\quad q^{\,7/8}\,\le\,X\,\le\,c_{2}q^{\,29/20}.
\end{cases}
\end{equation}
Since $3/4<97/128$ and $29/20<107/73$, then the estimates (\ref{lab_3-14}) and (\ref{lab_3-18}) give the bound
\begin{equation*}
\delta\,\ll\,
\begin{cases}
q^{\,6/7}X^{-8/7},& \text{if}\quad q^{\,3/4}\,\le\,X\,\le\,q^{\,7/8},\\
q^{\,16/35}X^{-24/35},& \text{if}\quad q^{\,7/8}\,\le\,X\,\le\,c_{1}q^{\,107/73}
\end{cases}
\end{equation*}
and
\begin{equation}\label{lab_3-19}
\Delta\,=\,\delta^{1/8}\,\ll\,
\begin{cases}
\bigl(q^{\,3/4}X^{-1}\bigr)^{1/7},& \text{if}\quad q^{\,3/4}\,\le\,X\,\le\,q^{\,7/8},\\
\bigl(q^{\,2/3}X^{-1}\bigr)^{3/35},& \text{if}\quad q^{\,7/8}\,\le\,X\,\le\,c_{1}q^{\,107/73}
\end{cases}
\end{equation}

It remains to consider the case
\begin{equation}\label{lab_3-20}
c_{1}q^{\,107/73}\,\le\,X\,\le\,(q/2)^{3/2}.
\end{equation}
Taking $V = 2Xq^{-1}$ we easily conclude that $1<V<\sqrt{X}$, $XV^{-1}=\tfrac{1}{2}\,q$, $V^{2}\le\tfrac{1}{2}\,q$ and
\[
\delta\,\ll\,\frac{D^{8}}{X^{8}}\,+\,\frac{\sqrt{q}}{D^{\,3/4}}\,+\,\frac{q^{2}}{DX}\,+\,\frac{X^{2}}{qD^{2\mathstrut}}.
\]
If we put $D = q^{\,2/35}X^{32/35}$ we find that
\[
\delta\,\ll\,q^{\,16/35}X^{-24/34}\,+\,q^{68/35}X^{-67/35}\,+\,q^{\,-39/35}X^{6/35}.
\]
After some calculations, we conclude that the first term dominates over two other terms for $X\gg q^{\,52/43}$ and for $X\gg q^{\,107/73}$ respectively.
Hence,
\begin{equation}\label{lab_3-21}
\delta\,\ll\,q^{\,16/35}X^{-24/34},\quad \Delta\,=\,\delta^{1/8}\,\ll\,\bigl(q^{\,2/3}X^{-1}\bigr)^{3/35}
\end{equation}
for any $X$ satisfying (\ref{lab_3-20}). The estimates (\ref{lab_3-19}) and (\ref{lab_3-21}) give the desired assertion. $\Box$
\vspace{0.3cm}

\textsc{Corollary.} \emph{Under the conditions of Theorem 1, we have}
\[
W_{q}(X)\,=\,W_{q}(a,b;X)\,=\,\mathop{{\sum}'}\limits_{p\le X}e_{q}(a\overline{p}+bp)\,\ll\,\pi(X)q^{\,2\vep}\Delta.
\]
\textsc{Proof.} By Abel summation, we find
\begin{multline*}
W_{q}(X)\,=\,\mathop{{\sum}'}\limits_{n\le X}\frac{\Lambda(n)}{\ln{n}}\,e_{q}(a\overline{n}+bn)\,-\,\sum\limits_{k\ge 2}\mathop{{\sum}'}\limits_{n = p^{k}\le X}\frac{\Lambda(n)}{\ln{n}}\,e_{q}(a\overline{n}+bn)\,=\\
=\,\mathop{{\sum}'}\limits_{n\le X}\frac{\Lambda(n)}{\ln{n}}\,e_{q}(a\overline{n}+bn)\,+\,O(\sqrt{X})\,=\,\frac{T_{q}(X)}{\ln{X}}\,+\,
\int_{2}^{X}\frac{T_{q}(u)\,du}{u(\ln{u})^{2}}\,+\,O(\sqrt{X}).
\end{multline*}
Estimating $T_{q}(u)$ trivially for $2\le u\le q^{\,3/4}$ and using the inequality $\sqrt{X}\ll q^{\,3/4}$ we obtain
\[
W_{q}(X)\,\ll\,\pi(X)q^{\,2\vep}\Delta\,+\,\frac{q^{\,3/4}}{(\ln{q})^{2}}\,+\,\int_{q^{\,3/4}}^{X}\frac{|T_{q}(u)|du}{u(\ln{u})^{2}}.
\]
If $X\le q^{\,7/8}$ then the last integral is estimated by
\[
\frac{1}{(\ln{X})^{2}}\int_{q^{\,3/4}}^{X}u^{\,-1/7}q^{\,3/28+2\vep}\,du\,\ll\,X^{6/7}q^{\,3/28+2\vep}(\ln{X})^{-2}\,\ll\,\pi(X)q^{\,2\vep}\Delta.
\]
If $q^{\,7/8}<X\le q^{\,3/2}$ then we split the segment of integration by the point $u = q^{\,7/8}$. Thus this integral is estimated as follows:
\begin{multline*}
(\ln{X})^{-2}\int_{q^{\,3/4}}^{q^{\,7/8}}u^{\,-1/7}q^{\,3/28+2\vep}\,du\,+\,(\ln{X})^{-2}\int_{q^{\,7/8}}^{X}u^{\,-3/35}q^{\,2/35+2\vep}\,du\,\\
\ll\,\bigl(q^{\,6/7}+X^{32/35}q^{\,2/35}\bigr)q^{\,2\vep}(\ln{X})^{-2}\,\ll\,X^{32/35}q^{\,2/35+2\vep}(\ln{X})^{-2}\,\ll\,\pi(X)q^{\,2\vep}\Delta.
\end{multline*}
It remains to note that
\[
\frac{q^{\,3/4}}{(\ln{q})^{2}}\,\ll\,\frac{X}{(\ln{X})^{2}}\,q^{\,3/4}X^{-1}\,\ll\,\frac{X\Delta^{7}}{(\ln{X})^{2}}\,\ll\,
\pi(X)\Delta^{7}.\quad \Box
\]

\pagebreak

\textbf{4. Applicatons.}

\vspace{0.3cm}

The above estimates allow one to establish the solvability of some congruences with prime numbers modulo $q$ lying in short interval $(1,N]$, $N\le q^{1-c}$, $c>0$.

\vspace{0.3cm}

\textsc{Theorem 2.} \emph{Let $0<\vep<0.01$ be any fixed number and suppose that $q\ge q_{0}(\vep)$ is prime, $(m,q)=1$. Further, let $q^{\,37/38+\vep}\le N\le \tfrac{1}{2}\,q$. Then the congruence}
\begin{equation}\label{lab_4-01}
p_{1}(p_{1}+p_{2}+p_{3})\,\equiv\,m\pmod{q}
\end{equation}
\emph{has a solution in primes $p_{1}, p_{2}, p_{3}$ such that $N<p_{1},p_{2},p_{3}\le 2N$.}
\vspace{0.3cm}

\textsc{Proof.} Since $(p_{1},q)=1$ then
\[
m\overline{p}_{1}\,-\,p_{1}-p_{2}-p_{3}\,\equiv\,0\pmod{q}.
\]
Setting $\pi_{1}(N)$ for the difference $\pi(2N)-\pi(N)$, and $I(N)$ for the number of solutions of (\ref{lab_4-01}),  we easily get
\begin{multline*}
I(N)\,=\,\frac{\pi_{1}^{3}(N)}{q}\,+\,R(N),\quad\text{where}\\
R(N)\,=\,\frac{1}{q}\sum\limits_{0<|c|<q/2}\biggl(\;\sum\limits_{N<p\le 2N}e_{q}(-cp)\biggr)^{\! 2}
\sum\limits_{N<p_{1}\le 2N}e_{q}(cm\overline{p}_{1}-cp_{1}).
\end{multline*}
Let $H$ be the maximal modulus of the sums
\[
\sum\limits_{N<p_{1}\le 2N}e_{q}(cm\overline{p}_{1}-cp_{1}),\quad 0<|c|<\tfrac{1}{2}\,q.
\]
Since $N>q^{\,37/38}>q^{\,7/8}$, the estimate of Theorem 1 implies that
$H\le Nq^{\,\vep/2}(q^{\,2/3}N^{-1})^{\,3/35}$ and hence
\[
|R(N)|\,\le\,\frac{Nq^{\,\vep/2}}{q}\,\bigl(q^{\,2/3}N^{-1}\bigr)^{\! 3/35}\sum\limits_{|c|\le q/2}\biggl|\;\sum\limits_{N<p\le 2N}e_{q}(-cp)\biggr|^{2}.
\]
The sum over $c$ is equal to $q\pi_{1}(N)$. Thus we get
\[
|R(N)|\,\le\,N^{2-3/35}q^{\,2/35+\vep/2}
\]
and therefore
\[
I(N)\,\ge\,\frac{\pi_{1}^{3}(N)}{q}\,\bigl(1\,-\,\delta(N)\bigr),
\]
where
\[
|\delta(N)|\,\le\,\frac{q}{\pi_{1}^{3}(N)}\,N^{2-3/35}q^{\,2/35+\vep/2}\,<\,N^{\,-38/35}q^{\,37/35+\vep}\,<\,q^{\,-3\vep/35}.
\]
Theorem is proved. $\Box$

\vspace{0.3cm}

\textsc{Theorem 3.} \emph{Let $k\ge 3$ and $0<\vep<0.01$ be any fixed constants and suppose that $q\ge q_{0}(\vep, k)$ is prime. Further, let $(ab,q)=1$ and $g(x)\equiv a\overline{x}+bx \pmod{q}$. Finally, suppose that}
\[
c_{k}\,=\,\frac{2k+31}{3k+29}\quad \textit{if}\quad 3\le k\le 9,\quad c_{k}\,=\,\frac{3k+22}{4(k+5)}
\quad \textit{if}\quad k\le 10
\]
\emph{and $q^{\,c_{k}+\vep}\le N\le \tfrac{1}{2}\,q$. Then the congruence}
\begin{equation}\label{lab_4-02}
g(p_{1})\,+\,\ldots\,+\,g(p_{k})\,\equiv\,m\pmod{q}
\end{equation}
\emph{has a solution in primes $p_{1},\ldots,p_{k}$ such that $p_{j}\le N$, $j = 1,\ldots, k$.}

\vspace{0.3cm}

\textsc{Proof.} Let
\[
H\,=\,\max_{(c,q)=1}\biggl|\,\sum\limits_{p\le N}e_{q}(cg(p))\biggr|,
\]
and let $I_{k}(N)$ be the number of solutions of (\ref{lab_4-02}). Then
\begin{multline*}
I_{k}(N)\,=\,\frac{\pi_{1}^{k}(N)}{q}\,+\,R_{k}(N),\\
|R_{k}(N)|\,\le\,\frac{1}{q}\sum\limits_{0<|c|<q/2}\biggl|\,\sum\limits_{p\le N}e_{q}(cg(p))\biggr|^{k}\,\le\,\frac{H^{k-2}}{q}\sum\limits_{|c|\le q/2}\biggl|\,\sum\limits_{p\le N}e_{q}(cg(p))\biggr|^{2}.
\end{multline*}
The sum over $c$ is equal to $q\kappa$, where $\kappa$ denotes the number of solutions of the congruence $g(p_{1})\equiv g(p_{2})\pmod{q}$ in primes $p_{1}, p_{2}\le N$. Since $q$ is prime, then the number of its solutions for fixed $p_{2}$ is at most two. Hence, $\kappa\le 2\pi(N)$ and $|R_{k}(N)|\le 2H^{k-2}\pi(N)$.

Now we note that $c_{k}\ge \tfrac{7}{8}$ for $3\le k\le 9$. Using the estimate of Theorem 1, for such $k$, we get the bound
\[
|R_{k}(N)|\,\le\,N^{k-1}(q^{\,2/3}N^{-1})^{\,3(k-2)/35}q^{\,\vep/2}\,=\,N^{\,(32k-29)/35}q^{\,2(k-2)/35+\vep/2},
\]
and therefore
\[
I_{k}(N)\,\ge\,\frac{\pi^{k}(N)}{q}\,\bigl(1\,-\,\delta_{k}(N)\bigr),
\]
where
\begin{multline*}
|\delta_{k}(N)|\,\le\,\frac{q}{\pi^{k}(N)}\,N^{\,(32k-29)/35}q^{\,2(k-2)/35+\vep/2}\,<\,q^{\,(2k+31)/35+\vep}N^{-(3k+29)/35}\,\le\\
\le\,q^{(2k+31)/35+\vep}\bigl(q^{\,c_{k}+\vep}\bigr)^{\,-(3k+29)/35}\,=\,q^{\,\vep-\vep(3k+29)/35}\,\le\\
\le\,q^{\,\vep-38\vep/35}\,=\,q^{-3\vep/35}.
\end{multline*}
Hence, $I_{k}(N)>0$.

If $k\ge 10$ then $\tfrac{3}{4}<c_{k}\le\tfrac{7}{8}$. Without loss of generality, we may assume that $N\le q^{\,7/8}$. Otherwise, the solvability of (\ref{lab_4-02}) in primes $p\le N$ follows from the solvability of (\ref{lab_4-02}) in primes $p\le M$, where $M=\bigl[q^{\,7/8}\bigr]$. By Theorem 1, we find that
\begin{multline*}
|R_{k}(N)|\,\le\,N^{k-1}\bigl(q^{\,3/4}N^{-1}\bigr)^{(k-2)/7}q^{\,\vep/2}\,=\,N^{(6k-5)/7}q^{\,3(k-2)/28+\vep/2},\\
I_{k}(N)\,\ge\,\frac{\pi^{k}(N)}{q}\,\bigl(1\,-\,\delta_{k}(N)\bigr),
\end{multline*}
where
\begin{multline*}
|\delta_{k}(N)|\,\le\,\frac{q}{\pi^{k}(N)}\,N^{(6k-5)/7}q^{\,3(k-2)/28+\vep/2}\,\le\,q^{\,(3k+22)/28+\vep}N^{-(k+5)/7}\,\le\\
\le\,q^{(3k+22)/28+\vep}\bigl(q^{c_{k}+\vep}\bigr)^{-(k+5)/7}\,=\,q^{\,\vep-\vep(k+5)/7}\,<\,q^{\,\vep-8\vep/7}\,=\,q^{\,-\vep/7}.
\end{multline*}
Theorem is proved. $\Box$
\vspace{0.3cm}

\textbf{5. Conclusion.}

\vspace{0.3cm}

In the above applications of Theorem 1, we consider only the prime modulus $q$. At the same time, Theorem 1 together with some additional estimates of ``long'' Kloosterman sums (with the length $X\gg q^{3/2}$) allow one to investigate the solvability of some congruences to any composite modulus. For example, we can prove
\vspace{0.3cm}

\textsc{Theorem 4.} \emph{Let $0<\vep<0.01$ be an arbitrary fixed constant and let $k\ge 3$ be any fixed integer. Suppose that $q\ge q_{0}(\vep, k)$. Further, let $(ab,q)=1$ and $g(x)\equiv a\overline{x}+bx \pmod{q}$. Finally, let}
\[
c_{k}\,=\,\frac{2(k+33)}{3k+64}\quad \textit{if}\quad 3\le k\le 16\quad\textit{and}\quad c_{k}\,=\,\frac{3k+50}{4(k+12)}
\quad \textit{if}\quad k\ge 17,
\]
\emph{and suppose that $q^{\,c_{k}+\vep}\le N\le q$. Then the number $I_{k}(N) = I_{k}(N,q,a,b,m)$ of solutions of} (\ref{lab_4-02}) \emph{in primes $p_{j}\le N$, $(p_{j},q)=1$, satisfies the relation}
\[
I_{k}(N)\,=\,\frac{\pi^{k}(N)}{q}\,\bigl(\varkappa_{k}(q)\,+\,O(\Delta_{k})\bigr).
\]
\emph{Here $\varkappa_{k}(q) = \varkappa_{k}(a,b,m;q)$ is some non-negative multiplicative function of $q$ for any fixed tuple $k,a,b$ and $m$. Moreover,}

a) \emph{for any $k\ge 7$ we have $\Delta_{k} = (\ln\ln{N})^{B}(\ln{N})^{-A}$,}
\[
A\,=\,\frac{1}{4}+\frac{57}{4}(k-7),\quad B\,=\,2^{k}-1;
\]

b) \emph{for any $k\ge 3$ we have $\Delta_{k} = q^{\,-\vep}$, if Generalized Riemann hypothesis is true.}
\vspace{0.3cm}

However, the structure of the ``singular series'' $\varkappa_{k}(a,b,m;q)$ (in the sense of dependence on the parameters $a,b, m$ and $q$) is quite complicated, especially for $k = 3,4$ and $q = 2^{n}$. We consider this question in a separate paper \cite{Changa_Korolev_2019}.

\renewcommand{\refname}{\normalsize{References}}

\end{document}